\documentclass[11pt,openany]{article}
\usepackage{amssymb}
\usepackage{amsmath}
\usepackage{amsthm}
\usepackage{amsfonts}
\usepackage{mathrsfs}
\usepackage{makeidx}
\usepackage{multicol}
\usepackage{fancyhdr}
\input{cyracc.def}

\newtheorem{thm}{Theorem}
\newtheorem{lemma}{Lemma}

\newtheorem{rem}{Remark}

\numberwithin{equation}{section}

\begin{document}

\title{\Large\bf  A note on a famous theorem of Pang and
Zalcman}

\author{Feng L\"{u}$^1$,\, Junfeng Xu$^2$\thanks{Corresponding author.}
\, and \,  Hongxun Yi$^3$\\
\small{1.College of Science,
China University of Petroleum,
Qingdao, 266580, P. R. China }\\
\small{ E-mail:  lvfeng18@gmail.com}\\
 \small{2. Department of Mathematics, Wuyi University, Jiangmen,
Guangdong 529020, P.R.China}\\
\small{ E-mail: xujunf@gmail.com}\\
\small{3. Department of Mathematics, Shandong University, Jinan,
250100 P.R. China.}\\
\small{E-mail:hxyi@sdu.edu.cn} }

\date{}
\maketitle
\thispagestyle{empty}

\begin{abstract}
In this paper, by studying the famous theorem of Pang and Zalcman,
we find a normal family and obtain a result, which is an improvement
of Pang and Zalcman's theorem in some sense. Meanwhile, several
examples are provided to show that our result's conditions are
necessary.

\end{abstract}

{\bf 2010 MSC:} Primary 30D45. Secondly 30D35.

{\bf Keywords and phrases:} normal family; share set; meromorphic
function;

\section{\bf Introduction}

Let $D$ be a domain in $\mathbb{C}$, let $f$ be a meromorphic
function on $D$, and let $S$ be a set with the finite elements. Set
$$
\overline{E}_f(S)=f^{-1}(\{S\})\cap D=\{z\in D:f(z)\in S\}.
$$

In this paper, we assume that $f,~g$ are two meromorphic functions
on $D$ and $S_1,~S_2$ are two sets. We denote
$\overline{E}_f(S_1)\subset \overline{E}_g(S_2)$ by $f(z)\in
S_1\Rightarrow g(z)\in S_2$. If
$\overline{E}_f(S_1)=\overline{E}_g(S_2)$, we denote this condition
by $f(z)\in S_1\Leftrightarrow g(z)\in S_2$. If the set $S$ has only
one element, say $a$, we denote $f(z)\in S$ by $f(z)=a$ (see
\cite{YY}).

Now, let $\mathcal{F}$ be a family of meromorphic functions on a
domain $D$. We say that $\mathcal{F}$ is normal in $D$ if every
sequence of functions $\{f_n\}\subset \mathcal{F}$ contains either a
subsequence which converges to a meromorphic function $f$ uniformly
on each compact subset of $D$ or a subsequence which converges to
$\infty$ uniformly on each compact subset of $D$(see. \cite{S}).

According to Bloch's principle, a lot of normality criteria have
been obtained by starting from Picard type theorems. On the other
hand, by Nevanlinna¡¯s famous five point theorem and Montel's
theorem, it is interesting to establish normality criteria by using
conditions known from a uniqueness theorem. A first attempt to this
was made by W. Schwick (see. \cite{Sc}).

Up to now, many normality criteria have been obtained in this
direction.(see. \cite{Ch,CF,CFZ,GM,LP,LY,LX,LXC,PZ,XZ,Xu}). In 2000,
Pang and Zalcman \cite{PZ} proved a famous theorem.\\

\noindent{ \sc\bf Theorem A}. {\it Let $\mathcal{F}$ be a family of
functions meromorphic on a domain, all of whose zeros are of
multiplicity (at least) $k$. If there exist $b\neq0$ and $h>0$ such
that for every $f\in \mathcal{F}$,
$\overline{E}_f(0)=\overline{E}_{f^{(k)}}(b)$ and
$0<|f^{(k+1)}(z)|\leq h$ whenever $z\in \overline{E}_f(0)$, then
$\mathcal{F}$ is normal in $D$.  }\\

It is natural to ask whether Theorem A still holds if the condition
$\overline{E}_f(0)=\overline{E}_{f^{(k)}}(b)$ is replaced by
$\overline{E}_f(0)\subset\overline{E}_{f^{(k)}}(b)$. Unfortunately,
we neither give a negative example nor prove it true. This problem
is very difficult even for the family of holomorphic functions(see.
\cite{Ch,CF,Xu}).  In this note, we study the special case that
$k=2$ and obtain the following result.

\begin{thm}\label{th1}
Let $\mathcal{F}$ be a family of functions holomorphic on a domain
$D$, all of whose zeros are of multiplicity (at least) $2$. If there
exist a non-zero constant $b$ and a positive constant $M$ such that
for every $f\in \mathcal{F}$,

(1) $f(z)=0\Rightarrow f''(z)=b$,

(2) $f''(z)=b\Rightarrow 0<|f'''(z)|\leq M$ and

(3) $f'^2(z)=Bf(z)$ whenever $z\in\overline{E}_{f''}(b)$,\\
where $B$ is a non-constant, then $\mathcal{F}$ is normal in $D$.
\end{thm}

\begin{rem}
Here, if $f$ omits a constant $b$, we can say that
all the zeros of $f-b$ are of multiplicity $\infty$.
\end{rem}

\begin{rem}
For the special cases that $\mathcal{F}$ is
holomorphic functions and $k=2$ of Theorem A, from
$\overline{E}_f(0)=\overline{E}_{f''}(b)$, it is easy to deduce
$\mathcal{F}$ satisfies the condition (3) of Theorem \ref{th1}.
Thus, in some sense, our result is an improvement of Theorem A.
Meanwhile, we know that the condition
$\overline{E}_f(0)=\overline{E}_{f^{(k)}}(b)$ is not necessary for
holomorphic functions in Theorem A.
\end{rem}

\begin{rem}
 We give an example to show that there exists a
normal family $\mathcal{F}$ satisfying the conditions of Theorem
\ref{th1}.
\end{rem}
Consider the family $\mathcal{F}=\{f_n,~n=1,~2,\ldots\}$ on the unit
disc, where
$$
f_n(z)=e^{\frac{z}{n}},
$$
so that
$$
f_n'(z)=\frac{1}{n}
~e^{\frac{z}{n}}~~~\hbox{and}~~~f_n''(z)=\frac{1}{n^2}
~e^{\frac{z}{n}}.
$$
Let $b$ be a non-zero constant and $B=b$. Then, it is easy to see
the family  $\mathcal{F}$ satisfies the conditions of Theorem
\ref{th1} and  $\mathcal{F}$ is normal on the unit disc.

\begin{rem}The assumption $0<|f''(z)|\leq M$ cannot be
replaced by $|f''(z)|\leq M$. We have a counter-example \cite{PZ} to
show it.
\end{rem}
 Consider the family
$\mathcal{F}=\{f_n,~n=1,~2,\ldots\}$ on the unit disc, where
$$
f_n(z)=\frac{1}{n^2}(e^{nz}+e^{-nz}-2)=\frac{1}{n^2}e^{-nz}(e^{nz}-1)^2,
$$
so that
$$
f_n^{(j)}(z)=n^{(j-2)}[e^{nz}+(-1)^je^{-nz}],~~j=1,~2,\ldots
$$
It is easy to see all the zeros of $f_n$ are of multiplicity 2 and
$$
f_n(z)=0\Leftrightarrow f''_n(z)=2\Rightarrow f_n'''(z)=0.
$$
While the family $\mathcal{F}$ is not normal on the unit disc.

\section{\bf Some Lemmas}

In order to prove our theorems, we need several lemmas.  For the
convenience of the reader, we recall these lemmas here.

The following result is due to Pang and Zalcman, see \cite{PZ}.
\begin{lemma}\label{le1}
Let $\mathcal{F}$ be a family of functions holomorphic on the unit
disc, all of whose zeros have multiplicity at least $k$, and suppose
that there exists $A \geq 1$ such that $|f^{(k)}(z)|\leq A$ whenever
$f(z)=0$, if $\mathcal{F}$ is not normal, then there exist, for each
$0 \leq \alpha \leq k,$

(a) a number $0 < r < 1;$

(b) points $z_{n},$ $z_{n}< r;$

(c) functions $f_{n}\in \mathcal{F} ,$ and

(d) positive number $\rho_{n}\rightarrow 0$ such that
$\rho_{n}^{-\alpha}f_{n}(z_{n}+\rho_{n}\xi)=g_{n}(\xi)\rightarrow
g(\xi) $ locally uniformly, where $g$ is a nonconstant holomorphic
function on $\mathbb{C}$, whose zeros have multiplicity at least
$k$,  such that $g^{\sharp}(\xi) \leq g^{\sharp}(0)=A+1$ and
$\rho(g)\leq 1$.

Here, as usual, $g^{\sharp}(\xi)=\frac{|g'(\xi)|}{1+|g(\xi)|^{2}}$
is the spherical derivative  and $\rho(g)$ is the order of $g$.
\end{lemma}
Next, we need to introduce a result, see \cite[Theorem 4.1]{HKR} or
\cite{NO}, which plays an important part in the proof of our
Theorem.
\begin{lemma}\label{le2}
Let $f$ be an entire function of order at most 1 and $k$ be a
positive integer, then
$$
m(r,\frac{f^{(k)}}{f})=o(\log~r),~~~\hbox{as}~~~r\rightarrow\infty.
$$
\end{lemma}
Finally, we recall the theorem of Chang, Fang and Zalcman, see
\cite{CFZ}, which is crucial to the proof of our theorem.
\begin{lemma}\label{le3}
Let $g$ be a non-constant entire function with $\rho(g)\leq 1$, let
$k\geq 2$ be an integer, and let $a$ be a non-zero finite value. If
$g(z)=0\Rightarrow g'(z)=a$, and $g'(z)=a\Rightarrow g^{(k)}(z)=0$,
then
$$
g(z)=a(z-z_{0}),
$$
where $z_{0}$ is a constant.
\end{lemma}

\section{\bf Proof of Theorem \ref{th1}}

Now, we prove Theorem \ref{th1}. For every $f\in \mathcal{F}$, it
follows from the assumption (1) that all the zeros of $f$ have
multiplicity 2. Noting that $f$ is holomorphic in $D$, we can set
\begin{equation}\label{3.1}
f=h^2,
\end{equation}
where $h$ is holomorphic in $D$. Differentiating $(\ref{3.1})$
yields
\begin{equation}\label{3.2}
f'=2hh',~~f''=2(h'^2+hh'')~~\hbox{and}~~f'''=6h'h''+2hh'''.
\end{equation}
We know that if $\mathcal{H}=\{h\}$ is normal in $D$, then
$\mathcal{F}$ is normal in $D$. Thus, we need only to prove that
$\mathcal{H}$ is normal in $D$. Suppose, to the contrary, that
$\mathcal{H}$ is not normal in $D$.

It is clear from (\ref{3.1}), the middle function of (\ref{3.2}) and
the condition (1) that
\begin{equation}\label{3.3}
h=0\Rightarrow h'\in \{a,~-a\}
\end{equation}
where $2a^2=b$. Combining the condition (2) and the last two
functions of (\ref{3.2}) yields
\begin{equation*}
2(h'^2+hh'')=b\Rightarrow 0<|6h'h''+2hh'''|\leq M.
\end{equation*}
By Lemma $\ref{le1}$, we can find $|z_{n}|<1$, $\rho_{n}\rightarrow
0$ and $h_{n}\in \mathcal{H}$ such that
\begin{equation}\label{3.5}
g_{n}(\xi)=\rho_n^{-1}h_{n}(z_{n}+\rho_{n}\xi)\rightarrow g(\xi)
\end{equation}
locally uniformly on $\mathbb{C}$, where $g$ is a non-constant
entire function such that $g^{\sharp}(\xi) \leq
g^{\sharp}(0)=M_1=|a|+1.$ In particular $\rho(g)\leq1$.

From $(\ref{3.5})$, it is easy to obtain that
\begin{equation}\label{3.6}
g_n'(\xi)=h_n'(z_n+\rho_n\xi)\rightarrow g'(\xi)
\end{equation}
and
\begin{equation*}
g_n''(\xi)=\rho_nh_n''(z_n+\rho_n\xi)\rightarrow g''(\xi)
\end{equation*}
locally uniformly on $\mathbb{C}$. Let
$$
H_n(\xi)=2[(g_n'(\xi))^2+g_n(\xi)g_n''(\xi)].
$$
Then, a routine calculation leads to
\begin{equation*}
H_n(\xi)=f''_n(z_n+\rho_n\xi).
\end{equation*}
Set
\begin{equation}\label{a}
G=2(g'^2+gg'').
\end{equation}
Thus, we can deduce that
\begin{equation}\label{3.7}
H_n(\xi)=2[(g_n'(\xi))^2+g_n(\xi)g_n''(\xi)]=f''_n(z_n+\rho_n\xi)\rightarrow
2[g'^2(\xi)+g(\xi)g''(\xi)]= G(\xi)
\end{equation}
locally uniformly on $\mathbb{C}$.

We claim that

(I)   ~~~$g(\xi)=0\Rightarrow g'(\xi)\in \{a,~-a\}$,

(II)  ~~~$g(\xi)=0\Rightarrow G(\xi)=b$ and

(III) ~~~$G(\xi)=b\Rightarrow G'(\xi)=0$.

First we prove (I).

Suppose that $g(\xi_{0})=0$, then by Hurwitz's theorem and
(\ref{3.5}), there exist a sequence $\{\xi_{n}\}$ such that
$\xi_{n}\rightarrow \xi_{0}$ and (for $n$ sufficiently large)
$$
g_{n}(\xi_{n})=\rho_n^{-1}h_{n}(z_{n}+\rho_{n}\xi_n)=0.
$$
Thus $h_{n}(z_{n}+\rho_{n}\xi_{n})=0$. It is clear from (\ref{3.3})
that
$$
h'_{n}(z_{n}+\rho_{n}\xi_{n})\in \{a,~-a\}.
$$
By (\ref{3.6}), we obtain
$$
g'(\xi_{0})=\lim\limits_{n\rightarrow\infty}h'_{n}(z_{n}+\rho_{n}\xi_{n})\in
\{a,~-a\},
$$
which implies $g(\xi)=0 \Rightarrow g'(\xi)\in \{a,~-a\}$. It is
(I).

Similarly as above, we can get (II).

We prove (III) as follows.

We affirm that $G\neq b$. Otherwise, suppose that $G=b$. That is
$$
2(g'^2+gg'')=b.
$$
Integrating the above differential equation yields $2gg'=bz+c$,
where $c$ is a constant.

If $g$ is a polynomial, then the equation $2gg'=bz+c$ implies that
$\deg(g)=1$. From (I), we get $g'= a$ or $-a $. Then
$$
|a|+1= g^{\sharp}(0)\leq |g'(0)|=|a| <|a|+1,
$$
a contradiction.

If $g$ is a transcendental entire function, then $g'$ is also a
transcendental entire function. By the lemma of logarithmic
derivative, we have
$$
\begin{aligned}
2T(r,g')&=T(r,g'^2)=m(r,g'^2)\leq
m(r,\frac{g'^2}{gg'})+m(r,gg')\\
&=m(r,\frac{g'}{g})+m(r,(bz+c)/2)=S(r,g)=S(r,g'),
\end{aligned}
$$
which is a contradiction. Thus, we finish the proof of $G\neq b$.

Now, we return to the proof of (III).

Suppose that $G(\zeta_{0})=b$. By Hurwitz's theorem and (\ref{3.7}),
there exist a sequence $\{\zeta_{n}\}$ such that
$\zeta_{n}\rightarrow \zeta_{0}$ and (for $n$ sufficiently large)
$$
H_{n}(\zeta_{n})=f''_n(z_n+\rho_n\zeta_n)=b.
$$
It follows from the assumption (2) that
$$
0<|f'''_n(z_n+\rho_n\xi_n)|\leq M.
$$
With (\ref{3.7}), we deduce
\begin{equation*}
H_n'(\xi)=\rho_nf'''_n(z_n+\rho_n\xi)\rightarrow G'(\xi)
\end{equation*}
locally uniformly on $\mathbb{C}$. Thus, it is not difficult to
deduce that
$$
G'(\zeta_0)=\lim\limits_{n\rightarrow\infty}
\rho_nf'''_n(z_n+\rho_n\zeta_n)=0,
$$
which implies (III).

Now, we continue to prove our theorem.

 Suppose that $\eta_0$ is a zero of $g$.
That is $g(\eta_0)=0$. By the claim (I) and (II), we get
$g'(\eta_0)=a~\hbox{or}~-a$ and $G(\eta_0)=b$. Differentiating
(\ref{a}) yields that
\begin{equation}\label{3.9}
G'=6g'g''+2gg'''.
\end{equation}
It is clear from (III) and (\ref{3.9}) that
$$
G'(\eta_0)=6g'(\eta_0)g''(\eta_0)+2g(\eta_0)g'''(\eta_0)=0.
$$
Then, we obtain $g''(\eta_0)=0$, which implies that
\begin{equation*}
g(\xi)=0\Rightarrow g''(\xi)=0.
\end{equation*}

Suppose that $g$ is a polynomial with $\deg{g}=n$. Noting that (I),
we know that $g$ has only simple zeros. Thus, $g$ has $n$ distinct
zeros $z_m$ ($m=1,~ 2,\ldots~,n$). By (I), we get $g'(z_m)=a$ or
$-a$ ($m=1,~ 2,\ldots~,n$). Thus, either $g'-a$ or $g'+a$ has at
least $p$ distinct zeros, here $p=\frac{n}{2}$ if $n$ is an even
number, $p=\frac{n+1}{2}$ if $n$ is an odd number. Without loss of
generality, we assume that $g'(z_m)-a=0$ ($m=1,~ 2,\ldots~,p$).
Obviously, $g''(z_m)=0$ ($m=1,~ 2,\ldots~,p$). It implies that each
$z_m$ ($m=1,~ 2,\ldots~,p$) is a multiple zero of $g'-a$.
Furthermore, it is easy to deduce that
$$
n-1=\deg(g')=\deg(g'-a)\geq 2p\geq n,
$$
a contradiction.

All the foregoing discussion shows that $g$ is a transcendental
entire function. Set
\begin{equation}\label{3.11}
\phi=\frac{g''}{g}.
\end{equation}
We find that $\phi$ is an entire function and
$\rho(\phi)\leq\rho(g)\leq 1$. Combining Lemma $\ref{le2}$ and the
lemma of logarithmic derivative yields
$$
T(r,\phi)=m(r,\phi)=m(r,\frac{g'}{g})=o(\log~r),
$$
which implies $\phi$ is a non-zero constant. By solving the
differential equation (\ref{3.11}), we have
\begin{equation}\label{3.12}
g=c_1e^{\lambda \xi}+c_2e^{-\lambda \xi},
\end{equation}
where $c_1,~c_2$ are two constants and $\lambda^2=\phi$.

Next, we prove that neither $c_1$ nor $c_2$ is zero. Otherwise,
without loss of generality, suppose that $c_2=0$. Combining
(\ref{a}) and (\ref{3.12}) yields
$$
G(\xi)=4c_1^2\lambda^2e^{2\lambda \xi}
$$
and
$$
G'(\xi)=8c_1^2\lambda^3e^{2\lambda \xi}.
$$
From (III) and the above two functions, it is easy to deduce a
contradiction. Thus, we finish the proof of that $c_1,~c_2$ are two
non-zero constants.

Differentiating the function $g$ yields
\begin{equation}\label{3.13}
g'(\xi)=\lambda [c_1e^{\lambda \xi}-c_2e^{-\lambda \xi}]
\end{equation}
and
\begin{equation}\label{3.14}
g''(\xi)=\lambda^2[c_1e^{\lambda \xi}+c_2e^{-\lambda \xi}].
\end{equation}
From (\ref{3.11}), it is obvious that
\begin{equation}\label{3.15}
g(\xi)=0\Leftrightarrow g''(\xi)=0.
\end{equation}
By (\ref{3.12}), we get
$$
g(\xi)=0\Leftrightarrow e^{\lambda\xi}\in\{A,~-A\},
$$
here $A= \sqrt{-\frac{c_2}{c_1}}$. From (I), we can see that
$$
e^{\lambda\xi}=A\Rightarrow g'(\xi)\in \{a,~-a\}.
$$
Noting that the form of $g'$, without loss of generality, we can
assume that
$$
e^{\lambda\xi}=A\Rightarrow g'(\xi)= a.
$$
Thus, we have
\begin{equation}\label{3.16}
g'(\xi)-a=e^{-\lambda\xi}[c_1\lambda
e^{2\lambda\xi}-ae^{\lambda\xi}-c_2\lambda]=A_1e^{-\lambda\xi}[e^{\lambda\xi}-A][e^{\lambda\xi}-A_2],
\end{equation}
where $A_1$ and $A_2$ are two non-zero constants. Observing that
(\ref{3.15}), we get
$$
e^{\lambda\xi}=A\Rightarrow g''(\xi)= 0,
$$
which implies that all the zeros of $e^{\lambda\xi}-A$ are multiple
zeros of $g'-a$. Therefore, we deduce that $A_2=A$. Rewriting
(\ref{3.16}) as
$$
g'(\xi)-a=A_1e^{-\lambda\xi}[e^{\lambda\xi}-A]^2.
$$
It indicates that $g'(\xi)=a\Leftrightarrow e^{\lambda\xi}=A$.
Meanwhile, with the same argument, we can deduce that
$g'(\xi)=-a\Leftrightarrow e^{\lambda\xi}=-A$. Combining the two
cases yields that $g'(\xi)\in \{a,~-a\}\Leftrightarrow
e^{\lambda\xi}\in \{A,~-A\}$. Thus, we have
$$
g(\xi)=0\Leftrightarrow g'(\xi)\in \{a,~-a\}.
$$
Furthermore, we obtain
\begin{equation}\label{3.17}
g=0\Leftrightarrow g'\in \{a,~-a\} \Leftrightarrow g''=0\Rightarrow
G=b.
\end{equation}
Noting that (\ref{3.13}), we know $g'-a$ has multiple zeros.
Differentiating $(\ref{3.14})$ yields
$$
g'''=\lambda^3 [c_1e^{\lambda \xi}-c_2e^{-\lambda \xi}].
$$
From the above function, it is not difficult to deduce that  $g'-a$
has zeros with multiplicity 2.

Suppose $g'(\alpha_0)=a$. By $(\ref{3.17})$ we get $g(\alpha_0)=0$
and $G(\alpha_0)=b$. From (III), we find that $\alpha_0$ is a
multiple zero of $G-b$. Noting that $G\neq b$, then there exists
$\delta>0$ such that
$$
g(\xi)\neq 0,~G(\xi)-b\neq0,
$$
in $D'(\alpha_0,\delta)=\{\xi:0<|\xi-\alpha_0|<\delta\}$. By
$(\ref{3.7})$, there exists $\varepsilon_0>0$ such that, for each
$0<\delta'<\delta$ and sufficiently large $n$,
$$
|f''_n(z_n+\rho_n\xi)-b-( G(\xi)-b)|< \varepsilon_0<|G(\xi)-b|
$$
on the circle $C(\alpha_0,\delta')=\{\xi:|\xi-\alpha_0|=\delta'\}$.
By Rouch$\acute{e}$ theorem, there exist $\{\alpha_{n,j}\}~(j=1,2)$
 tending to $\alpha_0$, such that, for each large $n$
\begin{equation}\label{3.18}
 H_{n}(\alpha_{n,j})=f''_n(z_n+\rho_n\alpha_{n,j})=b~~(j=1,2).
\end{equation}
And the assumption (2) implies that $\alpha_{n,1}\neq \alpha_{n,2}$.
Then, for $j=1,~2$, it follows from the assumption (3) that
\begin{equation}\label{A}
f_n'(z_n+\rho_n\alpha_{n,j})^2=Bf_n(z_n+\rho_n\alpha_{n,j}).
\end{equation}
We distinguish the
following three cases.

\textbf{Case 1.} For $j=1,~2$, there exist infinitely many $n_t^,$
satisfying
$$
f_{n_t}(z_{n_t}+\rho_{n_t}\alpha_{{n_t},j})=0.
$$
Then we get $h_{n_t}(z_{n_t}+\rho_{n_t}\alpha_{{n_t},j})=0$
$(j=1,~2)$. It follows from (\ref{3.5}) and Rouch$\acute{e}$ theorem
that $\alpha_0$ is a zero of $g$ with multiplicity at least $2$. But
$g$ has only simple zeros, a contradiction.

\textbf{Case 2.} For $j=1,~2$, there exist infinitely many $n_t^,$
satisfying
$$
f_{n_t}(z_{n_t}+\rho_{n_t}\alpha_{{n_t},j})\neq 0.
$$

We claim that there exists a subsequence of $\{n_t\}$ ( we still
denote it by  $\{n_t\}$) which contains infinite elements satisfying
\begin{equation}\label{B}
h_{n_t}'(z_{n_t}+\rho_{n_t}\alpha_{{n_t},j})=a ~~(j=1,~2).
\end{equation}
Without loss of generality, we need only to prove it holds for
$j=1$. By (\ref{3.1}), the first item of (\ref{3.2}) and (\ref{A}),
it is not difficult to deduce
$$
h_{n_t}'(z_{n_t}+\rho_{n_t}\alpha_{{n_t},1})\in \{d,~-d\},
$$
where $d=\frac{\sqrt{B}}{2}$ is a constant. It is clear from the
assumption $f_{n_t}(z_{n_t}+\rho_{n_t}\alpha_{{n_t},j})\neq 0$ that
$d$ is a non-zero constant.

Then, there must exists a subsequence of $\{n_t\}$ ( we still denote
it by
 $\{n_t\}$) which contains infinite elements satisfying
\begin{equation}\label{3.19}
h_{n_t}'(z_{n_t}+\rho_{n_t}\alpha_{{n_t},1})=e,
\end{equation}
here $e\in\{d,~-d\}$ is a non-zero constant. Then
$$
g'(\alpha_0)=\lim\limits_{n\rightarrow\infty}h_{n_t}'(z_{n_t}+\rho_{n_t}\alpha_{{n_t},1})=e.
$$
Noting that $g'(\alpha_0)=a$, we get $e=a$. With (\ref{3.19}), we
prove the claim.

On the other hand, by the middle item of (\ref{3.2}), (\ref{3.18}),
(\ref{B}) and the assumption of Case 2, we can deduce
$h_{n_t}''(z_{n_t}+\rho_{n_t}\alpha_{{n_t},j})=0$ for $j=1,2$.

Observing that $h_{n_t}''(z_{n_t}+\rho_{n_t}\alpha_{{n_t},j})=0$ for
$j=1,2$, so each $\alpha_{{n_t},j}$ $(j=1,2)$ is a multiple zero of
$h_{n_t}'(z_{n_t}+\rho_{n_t}\xi)-a$. It follows from (\ref{3.6}) and
Rouch$\acute{e}$ theorem that $\alpha_0$ is a zero
of $g'-a$ with multiplicity at least $4$, a contradiction.

\textbf{Case 3.} There exist  infinitely many $n_t^,$ satisfying
either
$$
f_{n_t}(z_{n_t}+\rho_{n_t}\alpha_{{n_t},1})=0,~~
f_{n_t}(z_{n_t}+\rho_{n_t}\alpha_{{n_t},2})\neq 0
$$
or
$$
f_{n_t}(z_{n_t}+\rho_{n_t}\alpha_{{n_t},1})\neq0,~~
f_{n_t}(z_{n_t}+\rho_{n_t}\alpha_{{n_t},2})=0.
$$

Without loss of generality, suppose that
$$
f_{n_t}(z_{n_t}+\rho_{n_t}\alpha_{{n_t},1})=0
~~\hbox{and}~~f_{n_t}(z_{n_t}+\rho_{n_t}\alpha_{{n_t},2})\neq 0.
$$
Similarly as Case 2, there exists a subsequence of $\{n_t\}$ ( we
still denote it by $\{n_t\}$) which contains infinite elements
satisfying
$$
h_{n_t}'(z_{n_t}+\rho_{n_t}\alpha_{{n_t},1})=a,
$$
$$
h_{n_t}'(z_{n_t}+\rho_{n_t}\alpha_{{n_t},2})=a~~\hbox{and}~~h_{n_t}''(z_{n_t}+\rho_{n_t}\alpha_{{n_t},2})=0.
$$
That means $\alpha_{n,2}$ is a multiple zero of
$h_{n_t}'(z_{n_t}+\rho_{n_t}\xi)-a$. Meanwhile,
$h_{n_t}'(z_{n_t}+\rho_{n_t}\xi)-a$ has another zero $\alpha_{n,1}$.
Then, it follows from (3.6) and Rouch$\acute{e}$ theorem that
$\alpha_0$ is a zero of $g'-a$ with multiplicity at least $3$, a
contradiction.

Thus, we get $g'(\alpha_0)\neq a$, which is a contradiction.

All the above discussion yields $\mathcal{H}$ is normal in $D$, so
$\mathcal{F}$ is also normal in $D$.

Hence, we complete the proof of Theorem \ref{th1}.

\section*{Acknowledgements}
This work was supported by  National Natural Science Foundation of China (Nos. 11126327, 11171184), NSF of Guangdong Province(No. S2011010000735) and the Outstanding Young Innovative Talents Fund of
Department of Education of Guangdong (No. 2012LYM0126).

\end{document}